\documentclass[11pt,reqno]{article}
\usepackage{amsmath, amssymb}
\numberwithin{equation}{section}

\newtheorem{theorem}{Theorem}

\newtheorem{conjecture}[theorem]{Conjecture}

\def \R {\mathbb{R}}

\def \s {\sigma}

\def \< {\langle}
\def \> {\rangle}

\def \vc {{\it v}}

\begin{document}
\title {Random processes via 
        the combinatorial dimension: introductory notes} 
\author {M. Rudelson\footnote{
   Department of Mathematics,
   University of Missouri,
   Columbia, MO 65211, USA;
   e-mail: \mbox{rudelson@math.missouri.edu} }
   \and
   R. Vershynin\footnote{
   Deptartment of Mathematics,
   University of California,
   Davis, CA 95616, USA;
   e-mail: \mbox{vershynin@math.ucdavis.edu}} }
\date{}

\maketitle

\begin{abstract}
  This is an informal discussion on one of the basic problems 
  in the theory of empirical processes, addressed in our preprint 
  {\em "Combinatorics of random processes and sections of convex bodies"},
  which is available at ArXiV and from our web pages. 
\end{abstract}

\qquad

Arguably, the central problem of the theory of empirical processes is
the following:
\begin{quote}
  {\em Describe the classes of functions on which 
  the classical limit theorems of probability hold uniformly.}
\end{quote}
Precisely, for a given probability space $(\Omega, \mu)$, one looks 
at the sequence of independent samples $(X_i)$ with values in $\Omega$ 
and distributed according to the law $\mu$. Then the problem is to describe 
the classes $F$ of functions $f : \Omega \to \R$ for which the sequence
of real valued random variables $f(X_i)$ satisfies the classical limit 
theorems of probability uniformly over $f \in F$. 
The two classical limit theorems we have in mind are 
the law of large numbers, for which such classes are called {\em Glivenko-Cantelli}, 
and the central limit theorem, which gives rise to {\em Donsker} classes.

In practice one does not know the law $\mu$ according to which the
samples $X_i$ are drawn. Thus a particularly vital question is to know what classes
are Glivenko-Cantelli or Donsker for every $\mu$. Such classes are called
{\em universal} or even {\em uniform} if the convergence in the 
corresponding limit theorems is uniform over all $\mu$.
We refer the reader to Chapter 14.3 of the book of Ledoux and Talagrand \cite{LT}
for a brief introduction to the theory of empirical processes and to the book
of Dudley \cite{Du 99} for a comprehensive account.

Vapnik and Chervonenkis were first to realize that the  
problem above is intimately connected with the combinatorics of the class $F$. 
In their pioneering works \cite{VC 68, VC 71, VC 81}, they quantified the 
size of a class $F$ by what now is called the {\em Vapnik-Chervonenkis dimension}.
This dimension makes sense only for Boolean classes, i.e. for classes of 
$\{0,1\}$-valued functions. A natural extension of this notion to general 
classes is possible (although not unique); the resulting concept is called 
by Talagrand ``the quantity of fundamental importance'' \cite{T 02}. 
This quantity measures how much $F$ oscillates. 

Precisely, for a given $t \ge 0$, a subset $\s$ of $\Omega$ 
is called $t$-shattered
by a class $F$ if there exists a level function $h$ on $\s$
such that, given any partition $\s = \s_- \cup \s_+$,
one can find a function $f \in F$ with $f(x)  \le  h(x)$ if
$x \in \s_-$ and $f(x)  \ge  h(x) + t$ if $x \in \s_+$.
The {\em combinatorial dimension} of $F$, denoted by $\vc(F,t)$,
is the maximal cardinality of a set $t$-shattered by $F$.
Simply speaking, $\vc(F,t)$ is the maximal size of a set
on which $F$ oscillates in all possible $\pm t/2$ ways
around some level $h$. 

For Boolean classes, this is clearly the classical 
Vapnik-Chervonenkis dimension in the whole nontrivial range $0 < t < 1$. 
For general classes, the combinatorial dimension appears 
implicitely in the work of Talagrand \cite{T 92} and 
explicitely in the paper of Alon et al \cite{ABCH}. 

Intuitively, the smaller the combinatorial dimension of 
a class $F$, the less $F$ oscillates, hence the better 
chances are for $F$ to be a Donsker class.
Recalling the known methods such as Dudley's integral inequality 
and Sudakov minoration, it is not difficult to locate the obstacles
for the class to be universal Donsker and, in fact, to guess the optimal 
description of such classes:   
\begin{conjecture}          \label{donsker}
  For every uniformly bounded class $F$, 
  $$
  \int_0^\infty \sqrt{\vc(F,t)} \; dt < \infty
  \ \Rightarrow \ \text{$F$ is universal Donsker}
  \ \Rightarrow \ \vc(F,t) = O(t^{-2}).
  $$
\end{conjecture}

The history of the work related to this conjecture goes back to 70's.
The right hand side, which makes sense for $t \to 0$, is known \cite{Du 99}.
The left hand side is difficult.  

In 1978, R.~Dudley proved Conjecture \ref{donsker} for classes of 
$\{0,1\}$-valued functions, where it clearly reads as follows: 
$F$ is a uniform Donsker class if and only if its combinatorial 
dimension $\vc(F,1)$ is finite (see \cite{Du 99}, \cite{LT} 14.3).
As Talagrand mentiones in \cite{LT} 14.3, this remarkable characterization 
is one of the main results on the combinatorics of empirical processes 
on $\{0,1\}$ classes.

In his 1992 Inventiones paper, Talagrand proved 
Conjecture \ref{donsker} for convex classes, however up to an additional 
factor of $\log^M(1/t)$ in the integrand  \cite{T 92}.
In 2002, Talagrand rewrote his proof 
in \cite{T 02}, now valid for arbitrary classes, and asked about the 
optimal value of the absolute constant exponent $M$. 

In their recent Inventiones paper \cite{MV}, Mendelson and Vershynin 
reduced the exponent of the logarithm in Talagrand's inequality 
to $M = 1/2$. This was achieved by a new argument, via 
construction of an unbalanced separated tree in $F$. 
However, these methods alone could not possibly remove the factor 
of $\sqrt{\log(1/t)}$. 

The logarithmic factor is removed completely in the forthcoming
preprint \cite{RV}. Thus the optimal exponent is $M = 0$ 
and Conjecture \ref{donsker} is true.
Proving this requires, in addition to the combinatorial ingredients of 
\cite{MV}, a new iteration argument and a general probabilistic reduction 
scheme. Along the proof of Conjecture~\ref{donsker}, its quantitative form 
would be especially helpful in practical applications. It can be written
in words as follows:
\begin{quote}
  {\em In Dudley's entropy integral inequality, the entropy can be 
  replaced by the combinatorial dimension.}
\end{quote}

The argument is based on the comparison of the combinatorial dimension and 
the uniform entropy $D(F,t)$ (also known as Koltchinskii-Pollard entropy), 
which is the supremum of the metric entropies of $F$ with respect to 
the metric of $L_2(\mu)$ over all probabilities $\mu$.
Even though the combinatorial dimension and the uniform entropy are 
not equivalent in general, the situation changes dramatically when we look
at {regular classes}, those for which, say, 
$\vc(F, 2t) \le \frac{1}{2} \vc(F,t)$ for all $t > 0$.  
Talagrand proved the required equivalence for these classes, however 
with an additional factor of $\log^M(1/t)$, where $M$ is an absolute constant
\cite{T 02}. In the Inventiones paper \cite{MV} of Mendelson and Vershynin, 
the exponent was reduced to $M = 1$ even {\em without} the regularity 
asumption on the class. In general, this exponent is optimal.

In \cite{RV} we prove that {\em with} the regularity assumption, the 
logarithmic factor can be removed completely; thus the optimal exponent 
in Talagrand's result is $M =0$, hence 
\begin{quote}
  {\em the uniform entropy and the combinatorial dimension are equivalent 
  for all regular classes.} 
\end{quote}
The argument is similar to that for Conjecture \ref{donsker}: 
the iteration scheme works in both cases, for the integral and for the 
regular quantity.

\end{document}